\documentclass[10pt]{article}
\usepackage{amssymb}
\usepackage{amsmath}

\begin{document}
\begin{center}
\Large
  Dembowski's Theorem on Finite Inversive Planes of Even Order
\end{center}

\begin{flushleft}
Bhaskar Bagchi\footnote{The author is a retired professor of Indian Statistical Institute.}  \\
B36 Medley, Good Earth, Malhar Avenue,\\
Kambipura, \\
Bangalore 560 074, \\
India.\\
\verb+bhaskarbagchi53@gmail.com+
\end{flushleft}

\begin{abstract}
A remarkable theorem due to Peter Dembowski states that if $I$ is an inversive plane of even order $q$ then $q$ must be a power of two and $I$ must be the incidence system of points versus plane ovals in an ovoid in the projective $3$-space over the field of order $q$. In this paper we present a short and self-contained proof of this result. Our proof depends on the classification  due to Benson of the symmetric and regular finite generalized quadrangles. Included here is a deduction of Benson's Theorem from the Dembowski-Wagner combinatorial characterization of finite projective geometries.
\end{abstract}

\section{Introduction}

A $t-(v,k, \lambda)$ design $D$ is an incidence system satisfying (i) $D$ has $v$ points, (ii) each block of $D$ is incident with exactly $k$ points, and (iii) any $t$ distinct points of $D$ are together incident with exactly $\lambda$ blocks of $D$. An easy counting argument shows that, for $0 \leq s \leq t$, any $t-(v,k,\lambda)$ design is also an $s-(v,k, \lambda_s)$ design where $\lambda_s = \lambda\binom{v-s}{t-s}/\binom{k-s}{t-s}$. (In the following, we shall use this formula several times, without further mention.) In particular, the number of blocks of the $t$-design is $b := \lambda_0 = \binom{v}{t}/\binom{k}{t}$, and, when $t>0$, the number of blocks incident with each point is $r :=\lambda_1= bk/v$.

A famous result of Fisher states that the parameters of any non-trivial $2$-design satisfy $b \geq v$ (equivalently $r\geq k$). A $2$-design with $b=v$ (equivalently $r=k$) is called a symmetric $2$-design. A $2$-design is symmetric iff any two distinct blocks of the design are together incident with $\lambda$ points. This happens iff the dual (obtained by interchanging the notion of points and blocks) of the $2$-design is again a $2$-design (necessarily with the same parameters). Note that a $2-(v,k,\lambda)$ design is symmetric iff its parameters satisfy $k(k-1)=\lambda(v-1)$.

The parameter $\lambda$ of a $t$-design is called its {\bf balance}. A $t$-design with balance $\lambda=1$ is said to be a {\bf  Steiner t-design}. A {\bf partial linear space} is an incidence system with at most one block incident with each pair of distinct points. In particular, when each point-pair is incident with a unique block, the incidence system is said to be a {\bf linear space}. The blocks of a partial linear space are usually called its {\bf lines}. Note that the Steiner $2$-designs are linear spaces. The symmetric Steiner $2$-designs are the {\bf finite projective planes}. These are precisely the $2$-designs with parameters $v=n^2+n+1, k = n+1, \lambda =1$. This number $n$ is called the {\bf order of the finite projective plane}. An {\bf affine plane of order n} is a $2-(n^2, n,1)$ design. Given a projective plane of order $n$ and a line $\ell$ in it, the incidence system obtained by deleting the line $\ell$ and the points on $\ell$ is an affine plane of order $n$. This process may be reversed as follows. It can be shown that the lines of an  affine plane naturally break up into $n+1$ "parallel classes" such that the lines in each parallel class partition the point set and any two lines from different parallel classes intersect (here we have adopted the usual convention of identifying any line of a partial linear space with the set of points incident with it). Given an affine plane of order $n$, one obtains a projective plane of order $n$ by adjoining $n+1$ points ("at infinity") corresponding to the parallel classes of the affine plane, and adjoining a single new line ("at infinity") incident with these points at infinity. This is called the {\bf projective closure} of the given affine plane.

Given an incidence system $D$, and a point $x$ of $D$, the {\bf contraction} $D_x$ of $D$ at $x$ is the incidence system whose points are the points of $D$ other than $x$, blocks are the blocks of $D$ incident with $x$, and whose incidence is the restriction of the incidence relation of $D$ to these points and blocks. Clearly, when $t>0$, each point contraction of a $t-(v,k,\lambda)$ design is a $(t-1)-(v-1, k-1, \lambda)$ design. An {\bf one-point extension} of a $t-(v,k, \lambda)$ design $D$ is a $(t+1)-(v+1, k+1, \lambda)$ $E$  (when it exists) such that $D$ is the contraction of $E$ at some point. An one-point extension of an affine plane of order $n$ is called an {\bf inversive plane of order n}. These are just the $3-(n^2+1, n+1, 1)$ designs. Clearly, the order of an inversive plane is the common order of all its point-contractions.
The blocks of an inversive plane are usually called the {\bf circles of the inversive plane}.

Given a prime power $q$, and $n \geq 1$. the n-dimensional projective space over the field of order $q$ is  denoted by $PG(n,q)$. We note that, when $n \geq 2$, the incidence system of points versus hyper-planes of $PG(n,q)$ is an example of a (symmetric) $2-(q^{n+1}-1)/(q-1), (q^n-1)/(q-1), (q^{n-1}-1)/(q-1))$ design. The projective space may be uniquely recovered from this design: the flats of the projective space are just the finite intersections of the blocks of the design. In particular, for prime powers $q$, $PG(2,q)$ is an example of a projective plane of order $q$.

An {\bf oval} in a projective plane of order n is a set of $n+1$ points no three of which are collinear. Thus each line of the plane meets an oval in $0,1\;\rm{or}\; 2$ points. Accordingly, the line is said to be a {\bf passant, tangent or secant line to the oval}. It is easy to see that each point of an oval is on a unique tangent line, so that there is a total of $n+1$ tangent lines. When $C$ is an oval in a projective plane of {\bf even order} n, there is a unique point $x$ of the plane such that the tangents to $C$ are precisely the lines through $x$. This point is called the {\bf nucleus} of the oval. Clearly, the nucleus of an oval does not belong to the oval.

 For a prime power $q$, an {\bf ovoid} in $PG(3,q)$ is a set of points such that every plane in the projective space meets this set in an oval or in a single point. A plane is said to be {\bf a tangent plane or a secant plane to the ovoid} according as it meets the ovoid in a point or in an oval It is easy to see that an ovoid in $PG(3,q)$ has exactly $q^2+1$ points in it, and every point $x$ in the ovoid is on a unique tangent plane to the ovoid. Indeed, the tangent planes to an ovoid constitute an ovoid in the dual projective $3$-space. It is easy to see that no three points in an ovoid are collinear. Since any three non-collinear points of $PG(3,q)$ are together in a unique plane, it follows that any three points of an ovoid are together in a unique oval contained in the ovoid.

 Given an ovoid $\mathcal{O}$ in $PG(3,q)$, let $I(\mathcal{O})$ denote the incidence system whose points are the points in $\mathcal{O}$, and blocks are the ovals contained in $\mathcal{O}$, It is immediate from the discussion in the previous paragraph that $I(\mathcal{O})$ is an example of an inversive plane (of prime power  order $q$). All the known finite inversive planes arise from this construction. The following is a famous theorem of Dembowski (see \cite{Dem1,Dem2}). Its proof occupies most of an entire chapter in Dembowski's book \cite{Dem2}, even though this book has an extremely cryptic style.

 {\bf Theorem 1.1} (Dembowski): Let $I$ be an inversive plane of even order $q$. Then $q$ is a power of two and $I$ is isomorphic to $I(\mathcal{O})$ for some ovoid $\mathcal{O}$ of $PG(3,q)$.

 It is widely believed that the order of any finite projective plane (equivalently, of any  finite affine plane) must be a prime power. Theorem 1.1 shows  that for an affine plane of even order to have an one point extension, the order must be a power of two. While it required powerful computers to establish that there is no affine plane of order $10$ (see \cite{Lam}), it is immediate from this theorem that there is no inversive plane of order $10$.

 Recall that a {\bf linear complex of lines} in $PG(3,q)$ is the set of all totally isotropic lines with respect to a non-degenerate symplectic form on (the underlying vector space of ) $PG(3,q)$. In the next section, we discuss regularity of generalized quadrangles and present a proof of Benson's theorem characterizing the linear complexes of lines in $PG(3,q)$ ($q$ prime power) as the only regular symmetric finite generalized quadrangles. Our proof depends on a beautiful combinatorial characterization of finite projective spaces due to Dembowski and Wagner. In the third and final section, we present a short proof of Dembowski's theorem (Theorem 1.1), deducing it from Benson's theorem. We believe that the proof given here is much more transparent than Dembowski's original proof.

 For unproved assertions on Design Theory/ Finite Projective Spaces made in this paper, the reader may consult \cite{CL,Dem2} and \cite{Hir}.

 \section{Benson's Theorem}

 Recall (\cite{PT}) that, for integers $s,t$, a {\bf generalized quadrangle} (in short, a GQ) of order $(s,t)$ is a partial linear space $X$ with $s+1$ points on each line and $t+1$ lines through each point satisfying the following property: given any point $x$ and line $\ell$ of $X$ such that $x \not \in \ell$, there is a unique point $y$ such that $y \in \ell$ and $y$ is collinear with $x$. An easy counting argument shows that a GQ $X$ of order $(s,t)$ has $(s+1)(st+1)$ points and $(t+1)(st+1)$ lines. The {\bf star of a point} $x$ (denoted $\rm{star}(x)$) in a GQ is the set of all points $y$ (including $x$ itself) collinear with $x$. Clearly the star of any point contains $s(t+1)+1$ points. It also readily follows from the definition of a GQ  that, for any two non-collinear points $x,y$ of $X$, exactly $t+1$ points are collinear with both $x$ and $y$. The set of these $t+1$ points is called {\bf the trace} of the pair $x,y$. Since any two of the points in the trace are non-collinear, it follows that at most $t+1$ points are collinear with all the points in the trace of $x,y$. This set of points is called {\bf the span} of the pair $x,y$. An unordered pair $\{ x,y \}$ of non-collinear points in a GQ $X$ is said to be  regular if the span of $\{ x,y \}$ has (the maximum possible) size $t+1$. $X$ is said to be a {\bf regular} GQ if all such pairs in $X$ are regular. Observe that if $\sigma$ is a span in a regular GQ, then the $t+1$ points in $\sigma$ are mutually non-collinear, and $\sigma$ is the span of any two of its points.

 The {\bf collinearity graph} of a partial linear space $X$  is  the graph having the points of $X$ as vertices, where two points are adjacent iff they are distinct and collinear in $X$. Recall that,for positive integers $m,n$, the {\bf complete bipartite graph} $K_{m,n}$ has $m+n$ vertices, split into two parts of size $m$ and $n$, such that two vertices of the graph are adjacent iff they belong to different parts, We have the following straightforward graphical reformulation of the notion of regularity of generalized quadrangles:

{\bf Lemma 2.1}:

The collinearity graph of any GQ $X$ of order $(s,t)$ has $ \leq\frac{1}{2} s^2(s+1)(st+1)(t+1)^{-1}$ induced subgraphs isomorphic to $K_{t+1,t+1}$. Equality holds here iff $X$ is regular.

{\bf Proof:} Let $N$ be the total number of induced $K_{t+1,t+1}$ in the collinearity graph of $X$. Note that, for any unordered pair $e$ of non-collinear points of $X$, and any induced subgraph $K_{t+1,t+1}$ containing $e$ in the collinearity graph of $X$, the two parts of this graph must be cthe trace and span of $e$. Thus $e$ is contained in at most one $K_{t+1,t+1}$, and such a $K_{t+1,t+1}$ exists iff $e$ is regular. Since $X$ has $(s+1)(st+1)$ points, and each point is non-collinear with $(s+1)(st+1)-(s(t+1)+1)= s^2t$ other points, it follows that there are $\frac{1}{2}s^2t(s+1)(st+1)$ such pairs $e$. Also each of the $N \; K_{t+1,t+1}$ in $X$ contains exactly
$2 \binom{t+1}{2}$ such pairs. Therefore we may count in two ways the total number of pairs $(e,G)$ where $e$ is as above and $G$ is a copy of $K_{t+1,t+1}$ whose vertex set contains $e$. This yields $2 \binom{t+1}{2}N \leq \frac{1}{2}s^2t(s+1)(st+1)$, with equality iff $X$ is regular. $\Box$

Let $q$ be a prime power. Recall that a {\bf polarity} of $PG(3,q)$ is an incidence preserving permutation interchanging points and planes (and hence mapping lines to lines). A point or plane is said to be absolute (with respect to a given polarity) if it is incident with its image under the polarity. A polarity is said to be a {\bf null polarity} if all points (equivalently planes) are absolute.  Let $W(q)$ denote the partial linear space whose points are the points of $PG(3,q)$, and whose lines are the lines of $PG(3,q)$ fixed by a given null polarity. Under its action by conjugation, the collineation group of $PG(3,q)$ is transitive  on the null polarities on $PG(3,q)$, so that this defines the incidence system $W(q)$ uniquely, up to isomorphism. A null polarity is given by ortho-complementation with respect to a non-degenerate symplectic form on the underlying vector space. Using this fact, it is easy to verify that, for any prime power $q$, $W(q)$ is a regular GQ of order $(q,q)$. In \cite{Ben}, Benson proved:

{\bf Theorem 2.2}  (Benson):

 If $W$ is a regular GQ of order $(q,q)$, then $q$ is a prime power and $W$ is (isomorphic to) $W(q)$.

 We recall that, if $x,y$ are two distinct lines of a $2$-design $D$, then {\bf the line of the 2-design} joining $x$ and $y$ is the intersection of all the blocks of $D$ containing $\{x,y\}$. If $\lambda$ is the balance of $D$, a line of $D$ may be defined as a set of $\geq 2$ points of $D$ which can be expressed as the intersection of $\lambda$ distinct blocks of $D$. A famous theorem of Dembowski and Wagner (\cite{DW}) states:

 {\bf Theorem 2.3}  (Dembowski and Wagner):

 Let $D$ be a symmetric 2-design of balance $>1$, Suppose every line of $D$ intersects every block of $D$. Then there is prime power $q$ and an integer $n \geq 3$ such that $D$ is isomorphic to the design of  points versus hyper-planes in $PG(n,q)$.

  {\bf Proof of Theorem 2.2}:  We begin with two claims. (i) every line of $W$ intersects every star in $W$, and (ii) every span in $W$ intersects every star in $W$. (i) is obvious since $W$ is a GQ: if $\ell$ is a line and $x$ is a point of $W$, $\ell \subseteq \rm{star}(x)$  when $x \in \ell$; and, when $x \not \in \ell$, $\ell \cap \rm{star}(x)=\{y\}$ where $y$ is the unique point in $\ell$ collinear with $x$. To prove (ii), note that, since $W$ is regular, all point pairs in a span $\sigma$ have a common trace, say $\tau$. If $x \in \tau$, then $\sigma \subseteq \rm{star}(x)$. in the contrary case, $\rm{star}(x)$ meets $\sigma$ in at most one point. Therefore, the sets $A(y) := \rm{star}(y)\setminus\tau, \; y \in \sigma,$ are $q+1$ pairwise disjoint sets of size $q^2$ each (If for $y_1 \not = y_2$ in $\sigma$, $z \in A(y_1)\cap A(y_2)$ then $z \not \in \tau$ and $z \in \rm{trace}(\{y_1,y_2\})= \tau$, contradiction.). They are contained in the complement of $\tau$ which is a set of size $(q+1)q^2$. So these sets partition the complement of $\tau$. Therefore, for any point $x \not \in \tau$, $x \in \rm{star}(y)$ (i.e., $y \in \rm{star}(x)$) for a unique point $y \in \sigma$. Hence $\sigma \cap \rm{star}(x) = \{y\}$. This proves (ii).

  Let $D$ be the incidence system whose points are the points of $W$ and whose blocks are the stars of points of $W$. Clearly $D$ has $(q+1)(q^2+1)$ points and equally many blocks. Also, as each point $x$ of $W$ is collinear with $1+q(q+1)=q^2+q+1$ points, it follows that each block of $D$ has size $q^2+q+1$. Also, for points $x_1 \not = x_2$, the blocks of $D$ containing both these points are $\rm{star}(x), x \in \rm{trace}(\{x_1, x_2\})$ when $x_1,x_2$ are non-collinear in $W$; when $x_1,x_2$ are collinear (say, lying on the line $\ell$ of $W$), these blocks are $\rm{star}(x), x \in \ell$. Thus  $D$ is a (symmetric) $2-((q+1)(q^2+1), q^2+q+1, q+1)$ design. Also, the line of $D$ joining any two collinear points is just a line of $W$, while the line of $D$ joining any two non-collinear points is the span of this pair. Therefore, in view of the preceding paragraph, every line of $D$ meets every block of $D$. Therefore, as $D$ is a symmetric design of balance $>1$, Theorem 2.3 applies and shows that $q$ is a prime power and $D$ is the points versus planes design of $PG(3,q)$.

  By this construction of $PG(3,q)$ from $W$, one sees that each plane of $PG(3,q)$ is of the form $\rm{star}(x)$ for a unique point $x$. But, for any two points $x,y$, we have $y \in \rm{star}(x)$ iff $x,y$ are collinear in $W$. Since collinearity is a reflexive and symmetric relation on the point set, it follows that  the map $x \mapsto {\rm star}(x)$ is a null polarity of $PG(3,q)$. From the description of the lines of $PG(3,q)$ given in terms of $W$ (as the lines and spans of $W$), it is now immediate that the lines of $W$ are the fixed lines of this polarity. Hence $W$ is $W(q)$. $\Box$

\section{Dembowski's Theorem}

An {\bf ovoid of a generalized quadrangle} is a set of points of the GQ meeting every line of the GQ in a unique point. Clearly an ovoid in a GQ of order $(s,t)$ has $st+1$ points, no two of which are collinear, and the star of a point $x$ meets it in $1$ or $t+1$ points -- according as $x$ does or does not belong to the ovoid.

{\bf Lemma 3.1}:

Any ovoid of $W(q)$ is an ovoid of the ambient $PG(3,q)$.

{\bf Proof:} Since the stars of points of $W(q)$ are just the planes in the ambient $PG(3,q)$, it follows that any plane of $PG(3,q)$ meets an ovoid $\mathcal{O}$ of $W(q)$ in one or $q+1$ points. Let $\ell$ be a line of $PG(3,q)$ meeting $\mathcal{O}$ in $m \geq 2$ points. The planes through $\ell$ induce a partion of the $(q^2+1-m)$-set $\mathcal{O} \setminus \ell$ into $q+1$ parts of size $q+1-m$ each. So we have $(q+1)(q+1-m)=q^2+1-m,$,i.e., $m=2$. Thus, no three points of $\mathcal{O}$ are collinear in $PG(3,q)$. It follows that $|mathcal{O}$ is an ovoid of $PG(3,q)$. $\Box$

Let $\mathcal{O}$ be an ovoid of $W(q)$. Since the secant planes to $\mathcal{O}$ are just the stars of the points in the complement of $\mathcal{O}$, we may identify the secant planes (and hence the ovals in $\mathcal{O}$) with the points of $W$ outside $\mathcal{O}$. Thus, the inversive plane $I(\mathcal{O})$ (defined in the introduction) has the following intrinsic description in terms of $W(q)$. The points and blocks of  $I(\mathcal{O})$ are the points in $\mathcal{O}$ and the points in the complement of $\mathcal{O}$, respectively. For points $x,y$ of $W(q)$ such that $x \in \mathcal{O}, y \not in \mathcal{O}$, $x$ is incident with $y$ in  $I(\mathcal{O})$ iff $x,y$ are collinear in $W(q)$.
 In view of Lemma 3.1, the following theorem is a version of Dembowski's Theorem (Theorem 1.1):

{\bf Theorem 3.2}:

Let $I$ be an inversive plane of even order $q$. Then $q$ is a power of two and $I$ is isomorphic to $I(\mathcal{O})$ for some ovoid $\mathcal{O}$ of $W(q)$.

Actually, when $q$ is a power of $2$ and $\mathcal{O}$ is an ovoid of $PG(3,q)$, it is easy to see that the partial linear space $W$, which has the points of $PG(3,q)$ as its points and the tangent lines to $\mathcal{O}$ as its lines, (is a regular GQ of order $(q,q)$, and hence) is isomorphic to $W(q)$. Clearly, the arbitrary ovoid of $PG(3,q)$, q even, is an ovoid of this copy of $W(q)$. (When $q$ is an odd prime power, $W(q)$ has no ovoids.) Thus, Theorem 3.1 is actually equivalent to Dembowski's Theorem.

The crucial observation behind the proof of Theorem 3.2 is the following. Since an inversive plane $I$ of order $q$ is a Steiner 3-design whose point-contractions are affine planes of order $q$, it follows that if $x$ is a point and $C$ is a circle (block) of $I$ such that $x \not \in C$, then $C$ is an oval in the projective closure $\pi$ of the contraction of $I$ at $x$. When $q$ is even, there is a point $y \not = x$ of $I$ (the nucleus of the oval $C$) through which all the tangent lines to $C$ in $\pi$  pass. (Since the line at infinity is a passant to $C$, the nucleus of $C$ must be an affine point.)

We say that two circles of an inversive plane are tangent if they have exactly one point in common. Recall (\cite{Dem2}, p. 253) that a {\bf pencil} in an inversive plane is a set of mutually tangent circles through some point $x$ which induces a partition of the set of all points other than $x$. The point $x$ is said to be the {\bf carrier} of the pencil.

{\bf Lemma 3.3}:

Let $I$ be an inversive plane of order $q$.

(a) If $x$ is a point and $C$ is a circle of $I$ such that $x \in C$, then there is a unique pencil $p$ of $I$ such that $C \in p$ and $x$ is the carrier of $p$.

(b) Suppose $q$ is even. If $p$ is a pencil of $I$ with carrier $x$ and $C$ is a circle of $I$ such that $ x \not \in C$ and $C \not \in p$ then there is a unique circle in $p$ which is tangent to $C$.

{\bf Proof}: (a) Note that a set $p$ of circles is a pencil of $I$ iff $ \{ C\setminus \{x\}: C \in p \}$ is a parallel class of lines in the affine plane $\pi$ obtained by contracting $I$ at $x$. Therefore, this result follows as every line of the affine plane is in a unique parallel class.

(b) Consider the contraction $\pi$ of $I$ at $x$. Let $\infty$ denote the point at infinity (in the projective closure of the affine plane $\pi$) corresponding to the parallel class $ \{D \setminus \{x\}: D \in p \}$. Note that $C$ is an oval in this plane of even order $q$. Since $\infty \not \in C$, and the line at infinity is not tangent to $C$ (obviously it is a passant to $C$) and passes through $\infty$, it follows that there is a unique line $l$ (in this projective closure) through $\infty$ which is a tangent to $C$. Then $(l \setminus \{\infty\})\cup \{x\}$ is the unique circle in the pencil $p$ which is tangent to $C$. $\Box$

{\bf Lemma 3.4}:

Let $e$ be a set of  two points in an inversive plane $I$ of even order $q$. Then there are exactly $q-1$ circles $C$ of $I$ such that $C$ is tangent to all the circles of $I$ containing $e$. These $q-1$ circles partition the complement of $e$ in the point set of $I$.

{\bf Proof}: Let $e=\{x,y\}$. Let $\mathcal{C}_e$ be the set of all circles $C$ of $I$ such that all the circles containing $e$ are tangent to $C$. We first note that, for $C \in \mathcal{C}_e$, we have $C \cap e = \emptyset$. Clearly, at most one of the points $x,y$ can be in $C$. Say $ y \not \in C$. Let $\pi$ be the affine plane obtained by contracting $I$ at $y$. Let $\overline{\pi}$ be its projective closure. $C$ is an oval in $\overline{\pi}$. Since $C \in \mathcal{C}_e$, all the lines of $\overline{\pi}$ through $x$ are tangent to $C$. Therefore, $x$ is the nucleus of $C$. Hence $x \not \in C$.

Now, let $C_1, C_2$ be two distinct circles in $\mathcal{C}_e$. We claim that $C_1 \cap C_2 = \emptyset$. Suppose not. Say, $z \in C_1 \cap C_2$.
Let $\pi$ be the contraction of $I$ at $z$. Let  $\overline{\pi}$ be the projective closure of $\pi$, obtained by adjoining the line at infinity $\ell$, say. Also, let $m$ be the line of  $\overline{\pi}$ joining the points $x,y$. For $i=1,2$, let $\ell_i$ be the line of  $\overline{\pi}$ extending the affine line $C_i \setminus \{z\}$. Let $w$ denote the point of  $\overline{\pi}$ at the intersection of $\ell_1, \ell_2$. Note that any of the $q$ circles $D$ containing $e$ but not containing $z$ is an oval of  $\overline{\pi}$. The point $w$ is in at least two tangents to $D$, namely $\ell_1, \ell_2$. Since  $\overline{\pi}$ is a projective plane of even order $q$, it follows that $w$ is the nucleus of $D$. So, $w \not \in D$ for any of the circles $D$ of $I$ such that $z \not \in D, e \subseteq D$. Also, as $\ell$ is a passant to $D$ and $m$ is a secant to $D$, it follows that $w \not \in \ell \cup m$. This is a contradiction since the lines $\ell,m$ together with these $q$ ovals $D$ clearly cover the entire point set of  $\overline{\pi}$. This proves the claim.

Thus, the circles in $\mathcal {C}_e$ are pairwise disjoint sets of size $q+1$ each, contained in the complement (of size $q^2-1$) of $e$. Therefore, $\# (\mathcal{C}_e) \leq (q^2-1)/(q+1)= q-1$, with equality iff the circles in $\mathcal{C}_e$ partition the complement of $e$. So, to complete the proof, it suffices to show that equality holds here for every 2-subset $e$ of the point set of $I$.

For any circle $C$ and point $x$ of $I$ such that $x \not \in C$, $C$ is an oval in the projective closure $\overline{\pi}$ of the contraction $\pi$ of $I$ at $x$. Since  $\overline{\pi}$ is of even order $q$, there is a unique point $y$ of  $\overline{\pi}$ such that all the lines of  $\overline{\pi}$ through $y$ are tangents to $C$. Since the line at infinity is a passant to $C$, $y$ must be an affine point. Thus, each of the $(q^2+1)-(q+1)=q^2-q$ points $x$ of $I$ in the complement of $C$ is in a unique 2-set $e=\{x,y\}$ such that $C \in \mathcal{C}_e$. Therefore, for each of the $q(q^2+1)$ circles $C$ of $I$, there are exactly $\frac{q^2-q}{2}$ 2-sets $e$ such that  $C \in \mathcal{C}_e$. Hence, as $e$ varies over the $\binom{q^2+1}{2}$ 2-subsets $e$ of the point set of $I$, the average size of $\mathcal{C}_e$ is $q(q^2+1)\frac{q^2-q}{2}/\binom{q^2+1}{2} = q-1$. But we have seen that $q-1$ is also an upper bound on the size of each $ \mathcal{C}_e$.
Hence $\#(\mathcal{C}_e)=q-1$ for all $e$. $\Box$

{\bf Proof of Theorem 3.2:}

Let $I$ be an inversive plane of even order $q$. Consider the incidence system $W$ whose points are the points and circles of $I$, and such that, corresponding to each pencil $p$ of $I$, $W$ has a block $\hat{p}$ consisting of the carrier of $p$ and the circles in $p$. From Lemma 3.3 (a), it is immediate that $W$ is a partial linear space with $q+1$ points on each line and $q+1$ lines through each point. Notice that, by the construction of $W$, we have (i) no two points of $I$ are collinear in $W$, (ii) if $x$ is a point and $C$ is a circle of $I$, then $x$ and $C$ are collinear in $W$ iff $x \in C$, and (iii) two circles $C_1, C_2$ of $I$ are collinear in $W$ iff $C_1$ and $C_2$ are (equal or) tangent circles. We claim that $W$ is isomorphic to $W(q)$. Once this claim is established, the proof will be complete since, by the construction of $W$, the point set $\mathcal{O}$ of $I$ is an ovoid of $W$ and $I=I(\mathcal{O})$.

Let $p$ be a pencil of $I$ with carrier $x$, and let $l=\hat{p}$ be the corresponding line of $W$. To prove that $W$ is a GQ (of order (q,q)), we need to show that any given point of $W$ not in the line $\l$ is collinear with a unique point of $W$ in $l$. Since the circles in the pencil $p$ induce a partition of the points other than $x$, this is obvious if the given point of $W$ is a point $y \not = x$ of $I$ : in this case  the unique circle in $p$ containing $y$ is the only point of $l$ collinear with $y$. Next let this point of $W$  be a circle  $C \not \in p$ of $I$. If $x \in C$, then $x$ is the only point on $l$ collinear with $C$ in $W$. If $x \not \in C$, then the circle of $I$ guaranteed by Lemma 3.3 (b) is the unique point in $l$ collinear with this point of $W$. Thus $W$ is indeed a GQ of order $(q,q)$.

With each 2-subset $e$ of the point set of $I$, we associate a subset of the point set of $W$ as follows. This set consists of the two points in $e$, the $q+1$ circles of $I$ containing $e$, and the $q-1$ circles of $I$ guaranteed by Lemma 3.4. It is immediate from Lemma 3.4 and the above description of collinearity in $W$ that the collinearity graph of $W$ induces a $K_{q+1,q+1}$ on this set. Thus, corresponding to the $\binom{q^2+1}{2}$ 2-subsets of the point set of $I$, we have found $\binom{q^2+1}{2}$ $K_{q+1,q+1}$ in the collinearity graph of $W$. But $\binom{q^2+1}{2}$ is the upper bound in Lemma 2.1 in the case $s=t=q$. Therefore, by Lemma 2.1, $W$ is a regular GQ of order $(q,q)$. Hence Theorem 2.2 implies that $q$ is a prime power and $W$ is isomorphic to $W(q)$. Since $q$ is even, it follows that $q$ is a power of two. $\Box$


\begin{thebibliography}{8}

 \bibitem{Ben} C.T. Benson, On the structure of generalized quadrangles,
 J. Algebra 15 (1970), 443-454.
 \bibitem{CL} P. J. Cameron and J.H. van Lint, Designs,Graphs,Codes and their Links,
 London Math. Soc. Student Texts 22, Cambridge Univ. Press, 1996.
 \bibitem{DW} P. Dembowski and A. Wagner, Some characterizations of finite projective spaces, Arch. Math. 11 (1960), 465-469.
 \bibitem{Dem1} P. Dembowski, Mobiusebenen gerader ordnung,
  Math. Ann. 157 (1964), 179-205.
  \bibitem{Dem2} P. Dembowski, Finite Geometries,
  Springer-Verlag,, New York, 1968.
  \bibitem{Hir} J.W.P. Hirschfeld, Finite Projective Spaces of Three Dimensions,
  Oxford Univ. Press, 1986.
  \bibitem{Lam} C.W.H. Lam, The search for a finite projective plane of order 10,
  Amer. Math. Monthly 98 (1991), 305-314.
  \bibitem{PT} S.E. Payne and J.A. Thas, Finite Generalized Quadrangles,
  Research Notes in Maths 110, Pitman, London, 1984.

\end{thebibliography}
 \end{document}